\documentclass{article}
\usepackage{amsmath,amsfonts,epsfig}
\usepackage{mathptmx}

\usepackage{amsmath,amsfonts,latexsym,amssymb}
\usepackage{mathrsfs}
\usepackage{verbatim}
\usepackage[utf8]{inputenc}
\usepackage[T1]{fontenc}
\usepackage{ae,aecompl}
\usepackage{braket}
\usepackage{color}
\newcommand{\bb}{\mathbb}

\usepackage{euscript}

\newcommand{\cx}{{\bb C}}

\newcommand{\hthree}{{\bb H}^3}
\newcommand{\htwo}{{\bb H}^2}

\newcommand{\ddz}{\frac{\partial}{\partial z}}

\renewcommand{\bold}[1]{\medskip \noindent {\bf \boldmath #1
                        }\nopagebreak[4]}

\newcommand{\del}{\partial}

\newcommand{\zbar}{{\overline{z}}}

\newcommand{\chat}{\widehat{\cx}}

\renewcommand{\Im}{\operatorname{Im}}

\renewcommand{\Re}{\operatorname{Re}}

\renewcommand{\skew}{\operatorname{skew}}

\newtheorem{theorem}{Theorem}[section]
\newtheorem{prop}[theorem]{Proposition}
\newtheorem{lemma}[theorem]{Lemma}
\newtheorem{cor}[theorem]{Corollary}

\newcommand{\cL}{{\cal L}}
\newcommand{\cM}{{\cal M}}

\newcommand{\calC}{{\mathcal C}}

\newtheorem{corollary}[theorem]{Corollary}
\newcommand{\ddt}{\frac{\del}{\del t}}

\renewcommand{\hbar}{\bar{{\mathbb H}}^3}

\newcommand{\CC}{\mathbb C}

\newcommand{\R}{\mathbb R}

\newcommand{\Half}{{\mathbb H}}
\newcommand{\Hp}{{\mathbb H}^2}
\newcommand{\Hs}{{{\mathbb H}^3}}

\newcommand{\barHs}{{\bar{\mathbb H}^3}}
\newcommand{\geu}{{g_{\rm euc}}}

\def\cotanh{\operatorname{cotanh}}

\def\eproof{$\Box$ \medskip}
\newcommand{\sltc}{\mathfrak{sl}_2(\mathbb C)}

\newcommand{\para}{{\mathfrak{p}}}
\newcommand{\pslt}{\mathsf{PSL}_2(\mathbb C)}

\newcommand{\Ep}{\operatorname{Ep}}
\renewcommand{\ddz}{\frac{\partial}{\partial z}}
\newcommand{\ddw}{\frac{\partial}{\partial w}}
\newcommand{\ddzbar}{\frac{\partial}{\partial \zbar}}
\newcommand{\E}{{\mathcal E}}
\newcommand{\Id}{{\operatorname{Id}}}
\newcommand{\ddx}{{\frac{\del}{\del x}}}
\newcommand{\ddy}{{\frac{\del}{\del y}}}

\begin{document}

\title{\bf \Large $L^2$-bounds for drilling short geodesics in  convex  co-compact hyperbolic 3-manifolds} \author{Martin
   Bridgeman\thanks{M. Bridgeman's research was supported by  NSF grants DMS-1500545, DMS-2005498 and institutional NSF grant DMS-1928930 while the author participated in the Fall 2020 program at MSRI. This work was also supported by a grant from the Simons Foundation (675497, Bridgeman)} \ and
  Kenneth Bromberg\thanks{K. Bromberg's research supported by NSF grant
     DMS-1509171, DMS-1906095.}}

\date{\today}

\maketitle
\begin{abstract}
We give $L^2$-bounds on the change in the complex projective structure on the boundary of conformally compact hyperbolic 3-manifold with incompressible boundary after drilling short geodesics. We show that the change is bounded by a universal constant times the square root of the length of the drilled geodesics. While $L^\infty$-bounds of this type where obtained in \cite{bromberg}, our bounds here do not depend on the injectivity radius of the boundary.  
\end{abstract}

\section{Introduction}
Given a complete, hyperbolic 3-manifold $M$ and a collection $\mathcal C$ of disjoint simple closed geodesics in $M$, the manifold $M - \mathcal C$ also supports a complete hyperbolic structure $\hat M$. If we insist that $M$ and $\hat M$ have the  same ending data then $\hat M$ is unique. If $M$ is closed, or more generally finite volume,  and the elements of $\mathcal C$ are sufficiently short then Hodgson-Kerckhoff \cite{Hodgson:Kerckhoff:rigidity} developed a theory of a hyperbolic cone-manifolds that allows one to continuously interpolate between $M$ and $\hat M$ through cone-manifolds. These methods were extended to conformally compact manifolds in \cite{bromberg:thesis}. By controlling the derivative of this family of cone-manifolds one can obtain comparisons between the geometry of $M$ and $\hat M$. 

One can give precise meaning to comparing the geometry of $M$ and $\hat M$. For example, one can compare the length of curves in $M$ to those in $\hat M$. In this paper  we will be interested in measuring the change in the projective boundary between the two manifolds. This change is described by a {\em holomorphic quadratic differential} given by taking the {\em Schwarzian derivative}. The size of this quadratic differential can be measured by taking an $L^p$-norm.  In \cite{bromberg}, the second author bounded the $L^\infty$-norm and these bounds played an important role the in resolving the Bers density conjecture. While $L^\infty$-bounds always imply $L^p$-bounds for all $p$, the bounds in \cite{bromberg} depended  on both the length of the curves being drilled and the injectivity radius of the hyperbolic metric on the boundary.  In this paper, we obtain $L^2$-bounds on the change in the projective structure  that are proportional to the square-root of the total length of the geodesics to be drilled but are independent of the injectivity radius. In particular, this gives uniform control on the $L^2$ change  for drilling short geodesics. In \cite{wp-paper}, this result is used to study the Weil-Petersson gradient flow of renormalized volume and  obtain lower bounds on the renormalized volume of a convex cocompact hyperbolic manifold  with incompressible boundary in terms of the Weil-Petersson  distance between its boundary components.
 
We have the following setup: $\bar N$ will be a compact, hyperbolizable 3-manifold with boundary with interior $N$ and $\mathcal C$ will be a collection of disjoint simple closed curves in $N$. The $\bar M$ is a complete, conformally compact hyperbolic structure on $\bar N$ where the curves in $\mathcal C$ are geodesics and $\bar M_t$ is a one parameter family of hyperbolic cone-manifolds with cone locus $\mathcal C$ and cone angles $t$. We also assume the conformal boundary $\del_c M_t$ is fixed throughout the definition.
\begin{theorem}[{\cite[Theorem 1.2]{bromberg}}]\label{cone_deformation}
There exists an $L_0>0$ such that if all geodesics in $\mathcal C$ have length $\le L_0$ in $M$ then the cone deformation exists for $t\in [0,2\pi]$ where $M_0$ is a complete, hyperbolic structure on $N - \mathcal C$.
\end{theorem}

While the conformal boundary will be a fixed conformal structure $X$ the deformation, the {\em complex projective structure} on $X$ will change. We denote this one parameter family of projective structures by $\Sigma_t$. The derivative of a path of projective structures on $X$ is naturally a holomorphic, quadratic differential. We denote the tangent vectors to $\Sigma_t$ by the holomorphic quadratic differentials $\Phi_t$. Our main results is the following bound on the $L^2$-norm of $\Phi_t$.

\begin{theorem}\label{inf_cone}
If $L_{\mathcal C}$ is the sum of the length the geodesics in $\mathcal C$ in $M = M_{2\pi}$ then
$$\|\Phi_t\|_2 \le c_{\rm drill} \sqrt{L_{\mathcal C}}.$$
\end{theorem}

As an immediate application we obtain the following $L^2$-bounds on the change in projective structure.

\begin{theorem}
There exists an $L_0>0$ and $c_{\rm drill}>0$ such that the following holds. Let $M$ be a conformally compact hyperbolic 3-manifold and $\mathcal C$ a collection of simple closed geodesics in $M$ each of length $\le L_0$. Let $\hat M$ be the unique complete hyperbolic structure on $M- \mathcal C$ such that the inclusion $\hat M \hookrightarrow \hat M$ is an isomorphism of conformal boundaries. If $\Sigma$ and $\hat\Sigma$ are the projective structures on the conformal boundaries of $M$ and $\hat M$ and the holomorphic quadratic differential $\Phi = \Phi(\Sigma, \hat\Sigma)$ is Schwarzian derivative between them then
$$\|\Phi\|_2 \le 2\pi c_{\rm drill} \sqrt{L_\mathcal C}$$
where $L_\mathcal C$ is the sum of the lengths of the components of $\mathcal C$ in $M$.
\label{cone}\end{theorem}
We note  that the $L^2$-bounds have universal constants compared to the $L^\infty$-bound in Theorem 1.3 in \cite{bromberg} which depended on injectivity radius of the boundary hyperbolic structure.  In \cite{bridgeman2019uniform}, $L^\infty$-bounds on quadratic differentials are obtained from $L^2$-bounds. These bounds again depend on the injectivity radius but they produce stronger bounds than those obtained in \cite{bromberg}. However, in \cite{bromberg} cone angles $> 2\pi$ where allowed which was important for the application to the Bers density conjecture.

\medskip

We briefly sketch our argument. Following the classical construction of Calabi \cite{calabi:rigidity} and Weil \cite{Weil:compact} the derivative of the deformation $M_t$ can be represented by a cohomology class in a certain flat bundle. This bundle has a metric and, in our setting, each cohomology class has a harmonic representative whose $L^2$-norm can be bounded by the length of the curves in the cone locus. We would like to use the bound on the $L^2$-norm in the 3-manifold to bound the $L^2$-norm of the quadratic differentials $\Phi_t$ representing the derivative of the projective structures.

To do this we first represent the cohomology class in the ends of the manifold by a certain {\em model deformation} which we describe explicitly. This model deformation will differ from the actual deformation by a trivial deformation. For the model deformation we can explicitly calculate the $L^2$-norm on the end in terms of the $L^2$-norm of the quadratic differential. To calculate the $L^2$-norm of the actual harmonic deformation we would like the model deformation to be orthogonal (in the $L^2$-inner product) to the trivial deformation. Unfortunately, this is not true, essentially because the model deformation itself is not harmonic. However, we will show that the model deformation is asymptotically harmonic. Using the infinitesimal inflexibility theorem from \cite{Brock:Bromberg:inflexible} we use this asymptotic control to bound the $L^2$-norm of the quadratic differential in terms of the $L^2$-norm of the deformation of the end.

This would seem to be enough however there is one final complication. Our bounds will depend on how large of an end we can embed in the cone-manifold. This size is controlled by the Schwarzian derivative of the of the projective structure. For smooth hyperbolic manifolds with incompressible boundary the Schwarzian of projective boundary is bounded by a classical theorem of Nehari. In the cone-manifold setting we will not be able to apply Nehari's theorem. Instead we control the Schwarzian by first controlling the {\em average bending} of the boundary of the convex core of the cone manifold. This notion was defined by the first author in \cite{Bridgeman:bending} where it was shown that for smooth hyperbolic 3-manifolds with incompressible boundary the average bending of the boundary of convex core is uniformly (and explicitly) bounded. We see that the argument in \cite{Bridgeman:bending} extends to cone-manifolds (with some restrictions) and then we will derive bounds on the Schwarzian via a compactness argument.

\subsubsection*{Acknowledgements} The work in this paper was motivated by a joint project with the authors and Jeff Brock. We thank Jeff for many interested discussions related to this paper.

\section{Background}
The proof relies on an analysis of $L^2$-bounds for cohomology classes associated to infinitesimal deformations of hyperbolic cone-manifolds. We now quickly review this theory with an emphasis and what is need for our computations. The original analysis can be found in \cite{bromberg:thesis} and \cite{bromberg} which generalized work of \cite{Hodgson:Kerckhoff:rigidity} on the finite volume hyperbolic cone-manifolds to the geometrically finite hyperbolic cone-manifolds.

Let $\barHs = \Hs\cup \chat$ be the usual compactification of $\Hs$ by $\chat$. Note that isometries of $\Hs$ extend to projective automorphisms of $\chat$ and that the group of isometries/projective transformations is $\pslt$. If $\bar M$ is a 3-manifold with boundary a $(\pslt,\barHs)$-structure on $\bar M$ is an atlas of charts to $\barHs$ with transition maps restrictions of elements of $\pslt$. On $M$, the interior of $\bar M$, this a {\em hyperbolic structure}. On the boundary $\del \bar M$ this is a {\em complex projective structure}.
In this paper we will be interested in a special class of $\left(\pslt,\barHs\right)$-structures, {\em conformally compact hyperbolic cone-manifolds}.

Let $\bar N$ be a compact 3-manifold with boundary with interior $N$ and let $\cal C$ be a collection of simple closed curves in the interior of $N$. Let $M = N - \calC$. A {\em hyperbolic cone metric} on $N$ with cone angle $\alpha$  along $\mathcal C$ is a hyperbolic metric on the interior of $M$ whose metric completion is homeomorphic to the interior of $N$ and in a c of each component the metric is that of a singular hyperbolic metric with cone angle $\alpha$. That is in cylindrical coordinates $(r, \theta, z)$ the metric will locally have the form
$$dr^2 + \sinh^2(r)d\theta^2 + \cosh^2(r)dz^2.$$
where $\theta$ is measured modulo the cone angle  $\alpha$ and the singular locus is identified with the $z$-axis. 

The hyperbolic metric is {\em conformally compact} if the hyperbolic structure on $M$ extends to a $\left(\pslt,\barHs\right)$-structure on $\bar M = \bar N - \cal C$. We then have:
\begin{theorem}[\cite{Hodgson:Kerckhoff:rigidity, bromberg:thesis}]
Given a cone angle $\alpha>0$ there exists a length $\ell>0$ such that the following holds. Let $(N,g)$ be a conformally compact hyperbolic cone-manifold with all cone angles $\alpha$ and assume that the tube radius about the singular locus is $\ge \sinh^{-1}\sqrt 2$. If each component of the singular locus has length $\le \ell$ then $t\in [0, \alpha]$ there a one parameter family of conformally compact hyperbolic cone-manifolds $(N, g_t)$ with the conformally boundary fixed and cone angle $t$.
\end{theorem}
This one parameter family of cone-manifolds will induces a one parameter family of projective structures $\Sigma_t$ on the boundary where the conformal structure of $\Sigma_t$ is fixed. We will be interested in controlling the change in this projective structure as the parameter varies.

\subsection{Flat $\sltc$-bundles}
The Lie algebra $\sltc$ can be interpreted geometrically as the space of infinitesimal automorphisms of $\barHs$. These are vector fields on $\barHs$ whose flow are elements in $\pslt$ so that on $\Hs$, the flow will be isometries of the hyperbolic metric, while on $\chat$ the flow will be projective automorphisms. A $\left(\pslt, \barHs\right)$-structure on $\bar M$ determines a flat $\sltc$-bundle $\bar E=E(\bar M)$ over $\bar M$. We examine this bundle when it is restricted to the hyperbolic structure and when it is restricted to the projective boundary.

\subsubsection*{Hyperbolic structures}
Let $M$ be the interior of $\bar M$. Then a $\left(\pslt, \barHs\right)$-structure is a hyperbolic structure and the bundle has a natural decomposition and metric structure that we now describe. Each fiber $E_p$ is the space of germs of infinitesimal isometries. In particular, $s \in E_p$ is a vector field in a neighborhood of $p$ so $s(p)$ is a vector in $T_p M$. As $E$ is a complex bundle we can multiply $s$ by $i$ and then $(is)(p)$ is another vector in $T_p M$. Then the map from $E$ to $TM\oplus TM$ given by $s\mapsto (s(p), (is)(p))$ is bundle isomorphism. In fact, the map $s\mapsto s(p) +i(is)(p)$ is a complex vector bundle isomorphism from $E$ to the the complexification $T^\CC M$ of the tangent bundle. This isomorphism from $E$ to $TM\oplus TM$ gives a decomposition of sections of $E$ into {\em real} and {\em imaginary} parts.

If $v$ is a vector in $T_p M$, we define $\hat v \in E_p$ such that under our isomorphism from $E_p$ to $T_p M \oplus T_p M$ we have $\hat v \mapsto (v, 0)$.  Then $\hat v$ is the infinitesimal translation with axis through $v$ and $i\hat v$ is an infinitesimal rotation about $v$. Note that $\hat v$ is real and $i\hat v$ is imaginary, as one would expect.

As $E$ is isomorphic to $TM\oplus TM$, the dual bundle $E^*$ is isomorphic to $T^*M \oplus T^*M$. The hyperbolic metric on $M$ determines an isomorphism from $TM$ to $T^*M$ and therefore an isomorphism from $E$ to $E^*$. Note that this isomorphism is $\R$-linear but is $\CC$-anti-linear with respect to the complex structures on $E$ and $E^*$. For sections $s$ of $E$ we let $s^\sharp$ be the dual section of $E^*$. When going from $E^*$ to $E$ we replace the $\sharp$ with a $\flat$.

As $E_p$ is a complex vector space  we have $E_p = E_p\otimes \CC$ and more generally for alternating tensors with values in $E_p$ we have
$$\Lambda^k(T_pM; E_p) = E_p \otimes \Lambda^k(T_p M) = E_p \otimes \Lambda^k\left(T^\CC_p M\right).$$
In particular, every $E$-valued form is locally the sum of terms $\phi s\omega$ where $\phi$ is a complex valued function, $s$ is a section of $E$, and $\omega$ is a $\R$-valued form. The $\sharp$ (and $\flat$) operators extend to $E$-valued forms and we have $(\phi s\omega)^\sharp = \bar\phi s^\sharp \omega$.
We also linearly extend the Hodge star operator from real forms to $E$-valued forms so that $\star(\phi s \omega) = \phi s(\star \omega)$.
 This extends to a linear map from $\Omega^k(M;E)$ to $\Omega^{3-k}(M;E)$. We the define the inner product
$$(\alpha, \beta) = \int_M \alpha \wedge (\star \beta)^\sharp.$$
Note that the wedge product of an $E$-valued form and an $E^*$-valued form is real form so this is a real inner product.  
We also let $\|\alpha\|^2 = (\alpha, \alpha)$ be the $L^2$-norm of an $E$-valued form.

If $\alpha$ is either an $E$-valued or $\CC$-valued form we define the pointwise norm $|\alpha|$ by $|\alpha|^2 = \star(\alpha \wedge (\star\alpha)^\sharp)$. Then $\|\alpha\|$ is the usual $L^2$-norm of the function $|\alpha|$.

If $d^*$ is the flat connection for $E^*$ then define the operator $\del$ on $E$ by
$$\del \omega = \left( d^*(\omega^\sharp) \right)^\beta.$$
Then the formal adjoint $\delta$ for $d$ satisfies the formula
$$\delta = \star \del \star.$$

 Note that if $s$ is a flat section ($ds =0$) then the real and imaginary parts, $\Re s$ and $\Im s$, will not be flat. That is $d$ will not preserve our bundle decomposition. Instead we define operators $D$ and $T$ such that $d = D+T$ where $D$ preserves the bundle decomposition and $T$ permutes it. That is for a real section $s$ we have that $Ds$ is a real $E$-valued 1-form while $Ts$ is imaginary. We have formulas for both $D$ and $T$. If $v$ is a vector field then 
$$D\hat v (w)= \widehat{\nabla_w v}$$
where $\nabla$ is the Riemannian connection for the hyperbolic metric and
$$T\hat v (w) =[\hat v, \hat w]$$
where $[,]$ is the Lie bracket. Note that the operator $T$ is purely algebraic.

We also have
$$\del = D-T.$$
This is a manifestation of the fact that the $\sharp$-operator is $\CC$-anti-linear.

The Laplacian for $E$-valued 1-forms is $\Delta = d\delta + \delta d$ and $\omega \in \Omega^k(E)$ is {\em harmonic} if $\Delta \omega = 0$. If $M$ is compact then this is equivalent $\omega$ be closed and co-closed. However, our manifolds will be non-compact so we will define $\omega$ to be a {\em Hodge form} if it is closed, co-closed and the real and imaginary parts are symmetric and traceless.

\subsubsection*{A computation in $E$}
We now make a few computations that will be very useful later and will also serve as an example of how to do computations in the bundle. 
We will work in the upper half space model of $\Hs = \CC\times \R^+$ with $\left\{\ddx, \ddy, \ddt\right\}$ and $\{dx, dy, dt\}$ the usual basis and dual basis at each $T_p\Hs$. We also let
$$\ddz = \frac12\left(\ddx - i\ddy\right) \mbox{ and } \ddzbar = \frac12\left(\ddx + i \ddy\right)$$
be tangent vectors in the complexified tangent space with dual 1-forms $dz = dx+i dy$ and $d\zbar = dx-idy$. We can then write any $E$-valued 1-form on $\Hs$ as a sum of $dz, d\zbar$ and $dt$ terms.

The Lie algebra $\sltc$ can be identified with traceless 2-by-2 matrices in $\CC$, projective vector fields on $\chat$ and infinitesimal isometries of $\Hs$. The reader can check the correspondence given in the following lemma:
\begin{lemma}\label{correspondence}
An element of $\sltc$ given by the matrix
$$\begin{bmatrix}
a& b\\
c & -a
\end{bmatrix}$$
is equivalent to the projective vector
$$2\left(-c z^2 + 2az + b\right)\ddz.$$
Along the axis $(0,t)$ in the upper half space model of $\Hs$
they are both equivalent to the constant section
$$2\left(ct^2\hat{\ddzbar} + a t\hat{\ddt}+ b\hat{\ddz}\right). $$
At $(0,t)$ we also have
$$\left|2\left(ct^2\hat{\ddzbar} + a t\hat{\ddt}+ b\hat{\ddz}\right)\right|^2 = 4|a|^2 +\frac{2|b|^2}{t^2}+2t^2|c|^2.$$
\end{lemma}

To calculate $Ts$ we note that
$$Ts = \left[s, \hat{\ddz}\right] dz + \left[s, \hat{\ddzbar}\right] d\zbar + \left[s, \hat{\ddt}\right] dt.$$

\begin{lemma}\label{Tparabolic}
Let $\para$ be a parabolic vector field on a neighborhood of a point $p$ in a hyperbolic manifold $M$. Let $e_n\in T_p M$ be a unit vector orthogonal to the horosphere tangent to $\para$, pointing away from the fixed point of $p$ and $\omega_n$ the dual $\R$-valued 1-form. Then
$$T\para = \hat e_n \otimes \omega + \para \otimes \omega_n$$
where $\omega$ is a $\CC$-linear 1-form with $\omega(e_n) =0$. Furthermore $|\omega| = |\para|$.
\end{lemma} 

{\bf Proof:} We can assume that $\para = \lambda z^2 \ddz$ and $p = (0,1)$ in the upper half space model. Then $e_n = \hat{\ddt}$ and $\omega_n = dt$. To calculate $T\para$ we use Lemma \ref{correspondence} to write $\para, \hat{\ddz}, \hat{\ddzbar}$ and $\hat{\ddt}$ as matrices and calculate the Lie brackets using matrix multiplication and get
$$T\para = \frac{\lambda}2\hat{\ddt} \otimes dz + \para\otimes dt$$
so $\omega = \lambda dz$. We can then compute to see that $|\omega| = |\para|$. \eproof

\subsubsection*{Complex projective structures}
We will be interested in complex projective structures $\Sigma$ that have a fixed underlying conformal structure $X$. The space $P(X)$ of such projective structures has a natural affine structure as the space $Q(X)$ of holomorphic quadratic differentials on $X$. That is the difference of $\Sigma_0$ and $\Sigma_1$ in is quadratic differential in $Q(X)$ defined as follows. Let $(U, \psi_0)$ and $(U, \psi_1)$ be charts for $\Sigma_0$ and $\Sigma_1$. Then $\psi_1\circ\psi_0^{-1}$ is a conformal map for an open neighborhood in $\chat$ to $\chat$. The {\em Schwarzian derivative} is a holomorphic function $\phi(\Sigma_0, \Sigma_1)$ on $\psi_0(U)$ and it determines a holomorphic quadratic differential $\Phi(\Sigma_0, \Sigma_1)$. Properties of the Schwarzian derivative imply that if $\Sigma_2 \in P(X)$ is a third projective structure then
$$\Phi(\Sigma_0, \Sigma_2) = \Phi(\Sigma_0, \Sigma_1) + \Phi(\Sigma_1, \Sigma_2).$$
This gives a canonical identification of the tangent space $T_\Sigma P(X)$ with $Q(X)$.

If $\hat g$ is a conformal metric on $X$ and $\Phi \in Q(X)$  then the ratio $|\Phi|/\hat g$ is a positive function on $S$. More concretely in a local chart $\Phi$ can be written as $\phi dz^2$ and the conformal metric can be written as $\hat g = \rho \geu$, where $\rho$ is a positive function and $\geu$ is the Euclidean metric. Then the ratio $|\phi|/\rho$, defined in the chart, is a well defined function on the surface. We denote this function by $\|\Phi(z)\|_{\hat g}$ and let $\|\Phi\|_{\hat g, p}$ be the $L^p$-norm of this function with respect to the $\hat g$ metric. We will mostly be interested in the hyperbolic metric but much of what we do will work in a more general setting. When we are using the hyperbolic metric we will drop the metric $\hat g$ from our notation.

For every conformal structure $X$ there is unique Fuchsian projective structure $\Sigma_F \in P(X)$. We let $\|\Sigma\|_{\hat g,\infty} = \|\Phi(\Sigma, \Sigma_F)\|_{\hat g,\infty}$.

If $\Sigma_t$ is a smooth path in $P(X)$ then its tangent vectors $\Phi_t$ lie in $Q(X)= T_{\Sigma_t}P(X)$. To prove our main result, Theorem \ref{cone}, we bound the distances between the endpoints of a path $\Sigma_t$ by bounding the norms of the derivative $\Phi_t$.

We also describe how a quadratic differential $\Phi \in Q(X) = T_\Sigma P(X)$ determines a cohomology class in $H^1(\Sigma; E(\Sigma))$. This was originally introduced by the second author in \cite{bromberg:thesis}.

Define a section $\para$ of $E(\CC)$ by
$$\para(z) = (w-z)^2\ddw = \frac12\begin{bmatrix} -z & z^2\\-1 & z\end{bmatrix}.$$
Then if $\Phi$ is represented in a projective chart by $\phi dz^2$ we defined an $E$-valued 1-form in the chart by
$$\para(z) \phi(z) dz.$$
One can then check that this gives a well defined $E$-valued 1-form $\omega_\Phi$ on $\Sigma$. As both $\para$ and $\phi$ are holomorphic, $\omega_\Phi$ is closed and therefore determines a cohomology class in $H^1(\Sigma; E(\Sigma))$.

\subsubsection*{Deformations of $\left(\pslt, \barHs\right)$-structures}
Let $v$ be a vector field on an open neighborhood $U$ in $\barHs$ that is conformal on $U\cap \chat$. We then define a section $s$ of $E(U) = U \times \sltc$ as follows.
For $x \in U\cap \Hs$ let $s(x)$ be the unique infinitesimal isometry that agrees with $v$ at $x$ and whose curl agrees with the curl of $v$ at $x$. On $U\cap \chat$, the vector field is (the real part) of the product of a holomorphic function $f$ and $\ddz$. For each $z\in U\cap \chat$ let $f_z$ be the complex quadratic polynomial whose 2-jet agrees with $f$ at $z$. Then $s(z) = f_z \ddz$.

Let $\left\{\left(U^\alpha, \psi^\alpha_t\right)\right\}$ be a 1-parameter family of $\left(\pslt, \barHs\right)$-structures on a 3-manifold with boundary $\bar M$ with the conformal boundary fixed. 
For each $x\in U^\alpha$, the time zero derivative of the path $\psi^\alpha_t(x)$ is vector field $v^\alpha$ on $\psi^\alpha_0(U^\alpha)$ that is conformal on $\psi^\alpha_0(U^\alpha)\cap \chat$. This determines a section $s^\alpha$ of $\bar E(\psi^\alpha_0(U^\alpha))$ and $\omega^\alpha=ds^\alpha$ is an $E$-valued 1-form on $\psi^\alpha_0(U^\alpha)$. While the sections $s^\alpha$ will not necessarily agree on overlapping charts, the $E$-valued 1-forms $\omega^\alpha$ will agree and determine an $E$-valued 1-form $\omega$ on $\bar M$. As locally $\omega$ is $d$ of a section, $\omega$ is closed and therefore represents an element of $H^1(\bar M; \bar E)$.

Since the conformal boundary is some fixed conformal structure $X$, the $\left(\pslt, \barHs\right)$-structures on $\bar M$ determine a family of projective structures $\Sigma_t \in P(X)$. The time zero derivative of $\Sigma_t$ will be a holomorphic quadratic differential $\Phi \in Q(X) = T_{\Sigma_0}P(X)$.
 \begin{prop}[{\cite[Theorem 2.3]{bromberg}}]
 The restriction of $\omega$ to the projective boundary is $\omega_\Phi$.
 \end{prop}

 \subsection{Hyperbolic metrics on ends}
 Let $\bar M$ have a $\left(\pslt, \barHs\right)$-structure. A convex surface $S$ in $M = \operatorname{int} \bar M$ {\em cuts off} a conformally compact end $\E$ if $\bar M - S$ has two components and the outward component is homeomorphic to $S \times [0,1)$ with $S\times \{0\} \subset \del \bar M$. Then $\E$ is the metric closure of the restriction of the outward component to $M$ and it is homeomorphic to $S\times (0,1]$ with $S = S\times\{1\}$ the original convex surface. We also let $\bar\E$ be the union of $\E$ with the projective boundary $\Sigma$ so that $\bar \E = S \times [0,1]$.
 
  The unit tangent vectors to the geodesic rays in $\E$ orthogonal to $S$ define a vector field on $\E$. We can choose this product structure such that $(z, s)$ is the time $s$ flow of this vector field. If $S_s$ is the the time $s$ image of $S$ under this normal flow, then the hyperbolic metric for $\E$ is can be written as a product of the induced metrics $g_s$ on $S_s$ and $ds^2$. 
 However, it will be convenient to parameterize the surfaces in the parameter $t = e^{-s}$ rather than in  the time parameter $s$ and we will see there is a nice formula for the hyperbolic metric in this product structure. The result is essentially due to Epstein. However, as it is not given in the exact form we need we derive it here.
 \begin{theorem}[C. Epstein, \cite{epstein-envelopes}]\label{model metric}
 Let $S$ be a convex surface cutting off a conformally compact end $\E$ with conformal boundary $X$. Then there exists a conformal metric $\hat g$ on $X$ and a bundle endomorphism $\hat B$ of $TX$ such that the hyperbolic  metric on $\E = S\times (0,1]$ is given by
 $$g_t \times dt^2/t^2$$
 where
 $$g_t = \frac1{4t^2} \left(\Id +t^2 \hat B\right)^* \hat g$$
 \end{theorem}
 
 {\bf Proof:} Let $g$ be the induced metric on $S$ and $B$ the shape operator and let $S_t$ be the distance $-\log t$ normal flow of $S$ in $\E$. We also let
 $$A_t = (1+t^2)\cdot \Id + (1-t^2)\cdot B = (\Id+B) +t^2\cdot (\Id-B).$$
Then the induced metric $g_t$ on $S_t$ is given by 
$$g_t = \frac1{4t^2} A^*_t g$$
 (see \cite[Lemma 2.2]{KS08}).  To
get our representation of the hyperbolic metric in $\E$ we need to rewrite $g_t$ in terms of a conformal metric on $X$.

The conformal structure on the boundary is induced from the conformal structure on $\E$. If we multiply our metric n $\E$ by $4t^2$ the conformal structure doesn't change but the new metric will extend continuously to 
$$\hat g = (A_0)^* g$$
 on  $S\times \{0\}$ so $\hat g$ is a conformal metric on $X$ and $A_0 = \Id+B$.
 As $S$ is convex, the eigenvalues of $B$ are non-negative. Therefore we define $\hat B = (\Id + B)^{-1}(\Id - B)$. It follows that 
$$g_t = \frac1{4t^2} A_t^* g = \frac1{4t^2} \left(A_0^{-1} A_t\right)^* \hat g = \frac1{4t^2} \left(\Id + t^2\hat B\right)^* \hat g$$
 as claimed. 
  \eproof
 
In Theorem \ref{model metric} the conformal metric on the boundary is determined by the convex surfaces in the hyperbolic end. In \cite{epstein-envelopes}, Epstein has a construction that starts with a metric at infinity and produces the convex surfaces. We only will use his construction for the hyperbolic metric. In particular, we have:
\begin{theorem}[{\cite[Propositions 6.4 and 6.5]{bromberg}}]
\label{epstein_hyp}
Let $X$ be the a component of the conformal boundary of a conformally compact hyperbolic cone manifold and $\hat g_x$ is the hyperbolic metric on $X$. Then for all $t> \log\left(\sqrt{1+ 2\|\Sigma\|_{\hat g_X,\infty}}\right)$ there is a convex surface $S$ that cuts of an end $\E$ such that $e^{2t}\hat g_X$ is the metric at infinity for $S$.
\end{theorem}

\subsubsection*{The metric in a chart}
If $(U, \psi)$ is a projective chart for $\Sigma$ then we can extend it to a chart $(U\times [0,1], \Psi)$  where $\Psi$ is the continuous extension of $\psi$ to a map to $\barHs$ that is an isometry on $U \times (0,1]$. We say that the chart is {\em adapted to $z_0 \in U$} if $\Psi(z_0,t)  = (0,t)$ where the coordinates on the right are in the upper half space model for $\Hs$. We can always construct a chart adapted to $z_0$ by taking any projective chart $(U,\psi)$ with $z_0 \in U$ and post-composing with an element of $\pslt$.

In a projective chart the metric at infinity $\hat g$ is scalar function times the Euclidean metric $\geu$. If the chart is adapted to $z_0$ then we can calculate the value of this function.

\begin{lemma}\label{adapted_chart}
If $(U\times [0,1], \Psi)$ is a chart adapted to $z_0 \in \Sigma$ then on $\psi(U)$ the metric $\hat g$ is of the form $\rho \geu$ where $\rho \colon \psi(U) \to \R^+$ is smooth and $\rho(0) = 4$.
\end{lemma}

{\bf Proof:} Define a function $\rho\colon \Psi(U\times [1,0)) \to \R^+$ with $\rho\circ\Psi(z,t) = 4t^2$.  If $g_\Hs$ is the metric for the upper half space model of $\Hs$ then $\rho\cdot g_\Hs$ extends continuously to the metric $\hat g$ on $\psi(U)$. Since $\Psi(z_0,t) = (0,t)$ we have that $\rho(0,t) = 4t^2$ and therefore $\rho g_{\Hs}$ extends continuously to $4\geu$ at $0$. \eproof

On a chart $(U, \psi)$ for $\Sigma$ we have the usual coordinate vector fields $\ddx$ and $\ddy$ along with the vector fields in $\ddz$ and $\ddzbar$ in the complexified tangent bundle. On a chart $(U\times [0,1], \Psi)$ for the end $\E$ these coordinate vector fields, along with $\ddt$ are a basis but, unlike in the upper half space model for $\Hs$, the vector fields $\ddx$ and $\ddy$ may not be orthogonal or of the same length as the operators $(\Id +t^2\hat B)$ are not conformal. In particular, the complex 1-form $dz$ will not be $\CC$-linear on the complex structure on $S$ induced by the metric $g_t$. However, we can write down $\CC$-linear and $\CC$-anti-linear forms in terms of the Beltrami differential of the endomorphisms that define $g_t$.

We begin with a computation on a single vector space. The usual Euclidean metric $\geu$ on $\R^2$ has a unique $\CC$-linear extension to $\R^2\otimes \CC$.
Then $dz$ and $d\zbar$ are the usual dual basis for $\R^2\otimes \CC$. While they are both $\CC$-linear on $\R^2\otimes \CC$, when restricted to $\R^2$, with the complex structure induced by $\geu$, $dz$ is $\CC$-linear while $d\zbar$ is $\CC$-anti-linear. A linear isomorphism $A\colon \R^2\to \R^2$ has a unique $\CC$-linear extension to $\R^2\otimes \CC$. The {\em Beltrami differential} for $A$ is
$$\mu = A_\zbar/A_z$$
where $A_z$ and $A_\zbar$ are complex numbers with
$$A^*(dz) = A_z dz + A_\zbar d\zbar.$$

We have the following:
\begin{lemma}{\label{hodge_star}}
Let $A\colon \R^2\to \R^2$ be a linear isomorphism and let $g = A^*\geu$ and $\mu$ the Beltrami differential for $A$. Then
$$\begin{bmatrix} dz\\ d\zbar\end{bmatrix} = \frac{1}{1-|\mu|^2} \begin{bmatrix} 1&-\mu\\ -\bar{\mu} & 1\end{bmatrix} \begin{bmatrix} dw\\d\bar{w} \end{bmatrix}$$
where  $dw$ and $d\bar w$ are $\CC$-linear and $\CC$-anti-linear on $\R^2$ with respect to $g$.
If $\star_g$ is the Hodge star operator for $g$ then
$$\star_g \begin{bmatrix}  dz \\ d\bar{z} \end{bmatrix} = \frac{-i}{1-|\mu|^2}\begin{bmatrix} 1 & \mu \\ -\bar{\mu}  & -1\end{bmatrix}\begin{bmatrix}  dw \\ d\bar{w} \end{bmatrix} = \frac{-i}{1-|\mu|^2}\begin{bmatrix} 1+|\mu|^2 & 2\mu \\ -2\bar{\mu}  & -1-|\mu|^2\end{bmatrix}\begin{bmatrix}  dz \\ d\bar{z} \end{bmatrix}.$$
Furthermore 
$$|dw|_{g} = \frac{|dz|_{g_{euc}}}{|A_z|}.$$
\end{lemma}

{\bf Proof:} As $dz$ and $d\zbar$ are $\CC$-linear and $\CC$-anti-linear for the complex structure on $\R^2$ induced by $\geu$, $A^*(dz)$ and $A^*(d\zbar)$ are $\CC$-linear and $\CC$-anti-linear for the complex structure induced by $g$. As $A$ is $\CC$-linear on $\R^2\otimes \CC$ and
$$A^*(dz) = A_zdz + A_{\zbar}d\zbar$$
for complex numbers $A_z$ and $A_\zbar$
we have
$$A^*(d\zbar) = \bar{A_{\zbar}} dz + \bar{A_z}d\zbar.$$
Dividing $A^*(dz)$ by $A_z$ and $A^*(d\zbar)$ by $\bar{A_z}$  we define $dw$ and $d\bar w$ by
$$\begin{bmatrix} dw\\d\bar{w} \end{bmatrix} = \begin{bmatrix} 1&\mu\\ \bar{\mu} & 1\end{bmatrix} \begin{bmatrix} dz\\ d\zbar\end{bmatrix}$$
so that $dw$ is $\CC$-linear and $d\bar w$ is $\CC$-anti-linear on the complex structure induced by $g$.
Inverting gives our formula for $dz$ and $d\zbar$ in terms of $dw$ and $d\bar w$.

As $dw$ and $d\bar w$ are $\CC$-linear and $\CC$-anti-linear with respect to $g$ we have
$$\star_g dw = -i dw, \qquad \star_g d\overline{w} = id\overline{w}$$
or
$$
\star_g \begin{bmatrix}  dw \\ d\bar{w} \end{bmatrix} = \begin{bmatrix} -i & 0\\ 0 & i\end{bmatrix}
\begin{bmatrix}  dw \\ d\bar{w} \end{bmatrix}.$$
We then have
$$\star_g \begin{bmatrix}  dz \\ d\bar{z} \end{bmatrix}
=  
\left(\frac{1}{1-|\mu|^2}\begin{bmatrix}  1& -\mu \\ -\bar{\mu} & 1\end{bmatrix}\right)
\begin{bmatrix} -i & 0\\ 0 & i\end{bmatrix}
\begin{bmatrix} 1 & \mu \\ \bar\mu & 1\end{bmatrix}\begin{bmatrix}  dz \\ d\bar{z} \end{bmatrix}
.$$
Multiplying, we obtain the stated formulas. 
For the norm $|dw|_{g}$ we note that define $dW = A^*dz = A_z dw$. Then
$$|A_z dw|_{g} = |dW|_{g} = |dz|_{g_{euc}}.$$
\eproof

We can apply the above to the metrics $g_t$. The Beltrami differential for endomorphisms $\frac1{2t}(\Id +t^2\hat B)$ can be written as $t^2\mu_t$ where 
$$\mu_t = \frac{\hat B_\zbar}{1 + t^2 \hat B_z}.$$ We obtain the following immediate corollary.

\begin{corollary}\label{model_beltrami}
Let $(U,\psi)$ be a projective chart for $\Sigma = \del \E$ with corresponding chart $(U\times[0,1], \Psi)$ for $\E$. Then
$$\star dz =-idz\wedge \frac{dt}{t} -2it^2\left(t^2\beta_0dz  +\beta_1 d\zbar\right)\wedge  \frac{dt}{t}.$$
where the $\beta_i$ are the smooth functions on $U\times [0,1]$ given by
$$\beta_0(z,t) = \frac{|\mu_t|^2}{1-t^4|\mu_t|^2} \qquad \beta_1(z,t) = \frac{\mu_t}{1-t^4|\mu_t|^2}.$$
 Further for $dw_t = dz+t^2\mu_t d\zbar$, 
$$\star dz = -idw_t\wedge  \frac{dt}{t} -it^2\left(t^2\beta_0dw_t  +\beta_1 d\bar w_t\right)\wedge  \frac{dt}{t}$$
and
$$|dw_t|_{g_t} = \frac{2t|dz|_{\hat g}}{|1+t^2B_z|} .$$
\end{corollary}

\subsection{Model deformations}
 If $S$ is a convex surface cutting of a conformally compact end $\E$ with projective boundary $\Sigma$ we can use Theorem \ref{model metric} to extend $\omega_\Phi$ to $\E$. Let
$$\Pi \colon S\times (0,1] \to \Sigma = S\times \{0\}$$
be given by $\Pi(z, t) = z$.
We would like to extend $\Pi$ to a bundle map between $E(\E)$ and $E(\Sigma)$. For this we note that for any flat bundle a path between two points in base determines an isomorphism between their fibers as the a flat bundle restricted to a path has a canonical product structure.

In our case the geodesic rays $\{z\} \times [0,1]$ are paths in $\bar{\E}$ between $(z,t)$ and $(z,0)$ and determine isomorphisms between the fiber $E_{(z,t)}$ of $E(\E)$ over $(z,t)$ and the fiber $E_z$ of $E(\Sigma)$ over $z$. 
Using this isomorphism we can extend $\Pi$ to a bundle map
$$\Pi_*\colon E(\E) \to E(\Sigma).$$
We then extend $\omega_\Phi$ to a 1-form in $\Omega^1(\E,E(\E))$ by pulling back $\omega_\Phi$ via $\Pi_*$.

 \begin{lemma}
 $$\omega_\Phi \wedge \star\left(\omega_\Phi^\sharp\right) = \frac{t^2}{16}\|\Phi\|^2_{\hat g}\left(1 + \frac{2t^4|\mu_t|^2}{1-t^4|\mu_t|^2}\right) dA_{\hat g}\wedge dt/t$$
 \end{lemma}
 	
 {\bf Proof:} We calculate at a point $(z,t) \in S\times (0,1]$ by taking a chart $(U, \psi)$ adapted to $z$. 
 In this chart $\omega_\Phi$ is written as
 $$\phi(z) \para(z) dz.$$
 While this expression does not depend on $t$, the Hodge star operator and the dual map will. In particular the expression
 $$\star({\phi}(z) \para(z) dz)^\sharp = \bar{\phi}(z) \para(z)^\sharp \star d\zbar$$
 depends on $t$ as both $p^\sharp$ and $\star d\zbar$ depend on $t$. 
 
 By taking the conjugate of $\star dz$ in Corollary \ref{model_beltrami} we have
 $$ dz \wedge \star d\zbar = i\left(1 + \frac{2t^4|\mu_t|^2}{1-t^4|\mu_t|^2}\right)dz \wedge d\zbar \wedge dt/t.$$

We also need
to find $\para(z)^\sharp(\para(z)) = |\para|^2(z)$. As we are working in a chart adapted to $z$, we have $\Psi(z,t)  = (0,t)$. Since $\para(0) = w^2 \ddw$, by Lemma \ref{Tparabolic} at $(0,t)$ we have
$$|\para(0)|^2 = \left|-t^2\hat\ddzbar\right|^2 = t^2/2.$$

By Lemma \ref{adapted_chart} we have $\|\Phi(z)\|^2_{\hat g} = 4|\phi(0)|^2$ and $dA_{\hat g} = 4dx\wedge dy = 2idz\wedge d\zbar$ and combining our calculations we have the result. \eproof

Let $\E_t = S \times (0,t]$ be portion of the end cutoff by $S_t$ and let $(,)_t$ be inner product on $\E_t$. For an $E$-valued form $\omega$ on $\E$ we then define $\|\omega\|^2_t = (\omega, \omega)_t$.
Integrating the prior lemma we immediately get:
\begin{cor}\label{L2_end_bounds}
We have $\|\omega_\Phi\|^2_t \ge 8t^2 \|\Phi\|^2_{\hat g, 2}$ and 
$$ \underset{t\to 0}{\lim} \frac1{t^2} \|\omega_\Phi\|^2_t = 8\|\Phi\|^2_{{\hat g},2}.$$
\end{cor}

The form $\omega_\Phi$ is not harmonic  as $\delta \omega_\Phi \neq 0$.  However, we will show that $\|\delta \omega_\Phi\|_2$ decays rapidly in $t$.

We'll break the estimate into small calculations.

\begin{lemma}\label{section_ext}
Let $U$ be a neighborhood of $0 \in \CC$
and let $s$ be a smooth section of $E(U \times [0,1])$ such that the function $|s|$ on $U\times (0,1] \subset \Hs$ extends continuously to $U \times [0,1] \subset \barHs$. Then the projective vector field $s(0,0)$ has a zero at $0 \in U$. If $|s|(0,0) = 0$ then $s(0,0)$ is a 
 is a multiple of $\para(0)$ and
$$|s|(0,t) = O(t).$$
\end{lemma}

{\bf Proof:} We write 
$$s(z,t) = \left(f_0(z,t) + f_1(z,t) w + f_2(z,t) w^2\right) \ddw$$
where the functions $f_i$ are smooth, complexed valued functions on $U\times [0,1]$. By Lemma \ref{correspondence} we have
$$|s|^2(0,t) = \frac{|f_0(0,t)|^2}{2t^2} + \frac{|f_1(0,t)|^2}{4} + \frac{t^2|f_2(0,t)|^2}{2}.$$
for $t \in (0,1)$. If $|s|$ extends continuously to $(0,0)$ then we must have $f_0(0,0) = 0$ so, as a projective vector field, $s(0,0)$ is zero at $0$.

If $|s|(0,0)  = 0$ the we must further have that $f_1(0,0) =0 $ and $\frac{\del f_0}{\del t} (0,0)=0$. Since $\frac{\del f_1}{\del t} (0,0)$ also exists it follows that 
$$|f_i(0,t)| = O\left(t^{2-i}\right)$$
and therefore
$$|s(0,t)| = O(t).$$
\eproof

\begin{lemma}\label{ds_ext}
Let $U$ be a neighborhood of $0 \in \CC$
and let $s$ be a smooth section of $E(U \times [0,1])$ such that the function $|s|$ on $U\times (0,1] \subset \Hs$ extends continuously to $U \times [0,1] \subset \barHs$. Then 
$$|ds \wedge d\zbar \wedge dt|(0,t) = O\left(t^3\right) \quad \mbox{and} \quad |ds \wedge dz \wedge dt|(0,t) = O\left(t^4\right).$$
\end{lemma}

{\bf Proof:} By Lemma \ref{section_ext} the condition that the norm $|s|$ extends to zero on $U\times \{0\}$ implies that $s(z,0) = f(z) \para(z)$ for some smooth complex valued function $f$ on $U$. We have that
$$ds\wedge d\zbar \wedge dt = s_z dz\wedge d\zbar \wedge dt = -2it^3 s_z dV$$
where $s_z$ is $z$-derivative of $s$ and $dV = \frac{idz\wedge d\zbar\wedge dt}{2t^3}$ is the volume form for $\Hs$. Therefore
$$|ds \wedge d\zbar \wedge dt|(0,t) = 2t^3|s_z|.$$
Note that the $z$-derivative of the section $\para$ is
$$\para_z(z) = -2(w-z)\ddw$$
and has bounded norm on $\Hs$. 
Therefore, on $\{0\} \times (0,1]$, the norm of $s_z = f_z\para + f\para_z$ is bounded and if we let $c$ be bound of $2|s_z|$ we have
$$|ds \wedge d\zbar \wedge dt|(0,t) \le c t^3.$$

The bound of the norm of $ds\wedge dz\wedge dt$ is similar once we note that the $\zbar$-derivative of $\para$ is zero so the norm of the $\zbar$-derivative $s_\zbar$  of $s$ is also zero at $(0.0)$. \eproof

We now prove our bounds on $\delta\omega_\Phi$.
\begin{lemma}
$$\lim_{t\to 0} \frac{1}{t^2}\|\delta\omega_\Phi\|_t = 0.$$
\label{delw_bound}
\end{lemma}

{\bf Proof:} We have $\delta = \star\del \star$ with $\del = D-T = d - 2T$. Therefore to bound $\|\delta \omega_\Phi\|_t$ we need to bound the norm of $d(\star\omega_\Phi)$ and $T(\star \omega_\Phi)$.

In a chart $(U \times [0,1], \Psi)$ by Corollary \ref{model_beltrami} we have
$$\star dz = -idz\wedge \frac{dt}{t} - 2it^2(t^2\beta_0 dz + \beta_1 d\zbar)\wedge \frac{dt}{t}$$
where the $\beta_i$ are smooth, complex valued functions on $U\times [0,1]$.
Therefore
$$\star \omega_\Phi = \phi \para \star dz = -i\phi\para dz\wedge \frac{dt}{t} -2it^2\phi\para(t^2\beta_0 dz + \beta_1 d\zbar)\wedge \frac{dt}{t}$$
As $\phi$ and $\para$ are holomorphic in the $dz$-coordinate we have that $d(-i\phi\para dz\wedge dt/t) = 0$. Therefore
$$d(\star \omega_\Phi)  = -2it(t^2d(\phi\beta_0\para) \wedge dz + d(\phi\beta_1\para)\wedge d\zbar) \wedge dt$$
and, as the sections $\phi\beta_i\para$ have norm limiting to zero on $\del \E$,  by Lemma \ref{ds_ext}
$$|d(\star \omega_\Phi)|(0,t) = O\left(t^4\right).$$

Next we calculate $T(\star\omega_\Phi)$. We will work in a chart adapted to $z$ and and a conformal coordinate $w_t$ at $(0,t)$. Again applying  Corollary \ref{model_beltrami} we have
$$\star\omega_\Phi =  \phi \para \star dz = -i\phi\para\left((1+t^4\beta_0) dw_t + t^2\beta_1 d\bar w_t\right)\wedge \frac{dt}{t}$$
where the $\beta_i$ are as above.
We use Lemma \ref{Tparabolic} to calculate $T\para$. At the point $(z,t)$ we have $e_n = \frac1t\ddt$, $\omega_n = \frac{dt}t$ and $\omega = \lambda dw_t$ for some scalar $\lambda$. Then 
$$T\para = \lambda t\hat\ddt\otimes dw_t + \para\otimes \frac{dt}{t}$$
where $|\lambda| |dw_t| = |\para| = t/\sqrt 2$. Then
$$T(\star\omega_\Phi) = -i\phi\lambda t^2\beta_1\hat\ddt\otimes dw_t \wedge d\bar w_t \wedge dt  = \frac{-i\phi t^4\beta_1  |dw_t|}{\sqrt 2} \hat\ddt dV$$
where $dV$ is the volume form and $dw_t \wedge d\bar w_t \wedge dt = it|dw_t|^2dV$. Then
$$|T(\star\omega_\Phi)| = t^3|\phi||\beta_1||dw_t|/\sqrt 2$$
since $\left|t\ddt\right| = 1$. By Corollary \ref{model_beltrami} $|dw_t| = O(t)$, giving
$$|T(\star\omega_\Phi)| =  O\left(t^4\right).$$
As the $\star$-operator is an isometry and $S$ is compact this implies that there exist a $c>0$ such that
$$|\delta\omega_\Phi|(z,t) \le ct^4.$$
We then have
\begin{eqnarray*}
\int_{\E_{t_0}} \delta\omega_\Phi \wedge \star\delta\omega_\Phi^\sharp &=& \int_{\E_{t_0}} |\delta\omega_\Phi|^2 dV \\
& \le & \int_0^{t_0}\int_S ct^8 dA_tdt/t\\
&\le & \int_0^{t_0} Kc t^5 dt = Kc{t_0}^6/6.
\end{eqnarray*}
Here $dA_t$ is the area form for the surface $S_t$ and we are using the fact that the area of these surfaces is bounded by $K/t^2$ for some $K>0$.
Therefore
$$\|\delta\omega_\Phi\|_t^2 \le Kct^6/6$$
and the lemma follows.
\eproof
  
 We now prove the main result of this section.

 \begin{theorem} 
 Let $\omega = \omega_\Phi +d\tau$ where the section $\tau$ of $E(\E)$ has finite $L^2$-norm.
The for all $t_0 \le 1$
$$||\Phi||^2_{\hat g,2} \leq \frac{1}{8t_0^2} ||\omega||_{t_0}^2$$
\label{hodge-limit}
\end{theorem}

{\bf Proof:}   We have 
$$(\omega,\omega)_t = (\omega_\Phi+d\tau, \omega_\Phi+d\tau)_t=  (\omega_\Phi, \omega_\Phi)_t+2\Re(d\tau,\omega_\Phi)_t+(d\tau, d\tau)_t.$$
By Corollary \ref{L2_end_bounds} we have
 $$\|\Phi\|_{\hat g,2}^2 \leq \frac{1}{8t^2} (\omega_\Phi,\omega_\Phi)_t.$$
 For the middle term we have integrate $\omega_\Phi\wedge \star d\tau^\sharp$ over the compact manifold $S\times [s,t]$ and let $s\to 0$. We have
\begin{eqnarray*}
\int_{S\times [s,t]} d\tau \wedge \star\omega_\Phi^\sharp &=& \int_{S\times [s,t]} \tau\wedge \delta\omega_\Phi + \int_{S_t} \tau\wedge\star\omega_\Phi^\sharp - \int_{S_s} \tau\wedge\star\omega_\Phi^\sharp\\
&=&\int_{S\times [s,t]} \tau\wedge \delta\omega_\Phi \to (\tau, \delta\omega_\phi)_t
\end{eqnarray*}
where the integrals over $S_t$ and $S_s$ are both zero since $\star\omega_\Phi$ restricted to these surfaces is zero as it contains a $dt$-term.

As $\|\tau\| < \infty$, applying Lemma \ref{delw_bound} we get
$$ \limsup_{t\rightarrow 0} \frac{1}{t^2}|(d\tau,\omega_\Phi)_t|  = \limsup_{t\rightarrow 0} \frac{1}{t^2}|(\tau,\delta\omega_\Phi)_t| \leq \limsup_{t\rightarrow 0} \frac{1}{t^2}\|\tau\|_t \cdot \|\delta\omega_\Phi\|_t= 0.$$
By the infinitesimal inflexibility theorem  \cite[Theorem 3.6]{Brock:Bromberg:inflexible} we have for $t < t_0$
$$ \frac{1}{t^2} (\omega,\omega)_t \leq \frac{1}{t_0^2} (\omega,\omega)_{t_0}.$$
Therefore as $(d\tau, d\tau)_t \geq 0$,
$$\|\Phi\|_{\hat g,2}^2 \leq \liminf_{t\rightarrow 0} \frac{1}{8t^2}(\omega_\Phi,\omega_\Phi)_t \leq \liminf_{t\rightarrow 0} \frac{1}{8t^2}(\omega,\omega)_t \leq  \frac{1}{8t_0^2} (\omega,\omega)_{t_0}$$
\eproof

If $\omega$ is a Hodge form on a conformally compact hyperbolic cone-manifold that is cohomologous to some $\omega_\Phi$ on and end  $\E$ then, by definition, $\omega= \omega_\Phi + d\tau$ for some $E$-valued section $\tau$ on $\E$. To apply this theorem we need the extra property that $\tau$ has finite $L^2$-norm. We call such a Hodge form a  {\em model Hodge form}.

\section{Nehari type bounds for cone-manifolds}
For a smooth, hyperbolic 3-manifold with incompressible boundary the classical Nehari bound on the Schwarzian derivative of univalent maps gives that $\|\Sigma\|_\infty \le 3/2$ for every component $\Sigma$ of the projective boundary. We are interested in obtaining similar bounds for a hyperbolic cone-manifolds. To do so we need to make some technical assumptions, that will always be satisfied in our applications, but do make the statement somewhat cumbersome.

One of the difficulties is that the usual Margulis lemma does not hold for cone-manifolds. The following statement is a replacement.
\begin{theorem}[{\cite[Theorem 3.5]{bromberg}}]\label{cone_margulis}
There exists an $L_0>0$ such that the following holds. Let $M$ be a hyperbolic cone-manifold such that all cone angle $\le 2\pi$, every component of the cone cone locus has length $\le L_0$ and  every component has a tubular neighborhood of radius $\sinh^{-1}\sqrt{2}$. Further assume that these neighborhood are mutually disjoint. Then each component $c$ of the cone locus of length $L_c$ and cone angle $\theta_c$ has a tubular neighborhood of radius $R_c$ where
$$\theta_c L_c \sinh(2 R_c) = 1.$$
\end{theorem}

Now we state our version of the Nehari bound. When the cone angle is small it will be important that the cone locus has a large tubular neighborhood where the radius grows as the cone angle decreases. The necessary lower bounds will come from the previous result and to use it we will need to assume that the length of the cone locus is bounded above by a linear function of the cone angle.
\begin{theorem}
There exists an $L_0>0$ such that the following holds. Let $M$ be a conformally compact hyperbolic cone-manifold such that all cone angle $\le 2\pi$ and there are a disjoint collection of tubular neighborhoods of the components of the cone locus of radius $\ge \sinh^{-1} \sqrt 2$. Further assume that if $c$ is a component of the cone locus with cone angle $\theta_c$ and length $L_c$ then
$$L_c \le \theta_c L_0.$$
 Then for every component $\Sigma$ of the projective boundary of $M$ we have
$$\|\Sigma\|_{\infty} \le K.$$
\label{infinity_bound}
\end{theorem}

In order to prove this, we will need to consider the Thurston parametrization of projective structures via measured laminations and use the notion of average bending of a measured lamination. We show that the result follows from a compactness argument.
\subsection{The Thurston parameterization}

The space $P(\Delta)$ of projective structures on the hyperbolic disk is equivalent to the space of locally univalent maps $f\colon\Delta\to \chat$ with the equivalence  $f\sim g$ if $f =\phi\circ g$ for some $\phi \in \pslt$. We can identify $P(\Delta)$ with the space of quadratic differentials $Q(\Delta)$ by mapping $[f] \in P(\Delta)$ to its Schwarzian derivative $S(f) \in Q(\Delta)$. Then the topology on $P(\Delta)$ is the compact-open topology on $Q(\Delta)$.

 Thurston described a natural parameterization of $P(\Delta)$ by $\cM\cL(\htwo)$ the space of measure geodesic laminations on $\htwo$. We briefly review this construction.

A round disk $D \subset \chat$ shares a boundary with a hyperbolic plane $\htwo_D \subseteq \hthree$. Let $r_D\colon D\to\hthree$ be the nearest point projection to $\htwo_D$ and $\tilde{r}_D\colon D\to T^1\hthree$ be the normal vector to $\htwo_D$ at $r_D(z)$ pointing towards $D$. We can 
use these maps to define a version of the Epstein map for $\rho_f$. In particular define $\widetilde\Ep_{\rho_f}\colon\Delta\to T^1\hthree$ by $\widetilde\Ep_{\rho_f}(z) = \tilde r_{f(D)}(f(z))$ where $D$ is the unique round disk with respect to $f$ such that $\rho_D(z) = \rho_f(z)$ and let $\Ep_{\rho_f}(z) = \pi\circ\widetilde\Ep_{\rho_f}(z) = r_{f(D)}(f(z))$. (For the existence of this disk see \cite[Theorem 1.2.7]{KT-projective}.) We also define $\widetilde\Ep_{e^t\rho_f} = g_t\circ\widetilde\Ep_{\rho_f}$ and $\Ep_{e^t\rho_f} = \pi\circ\widetilde\Ep_{e^t\rho_f}$.

The image of $\Ep_{\rho_f}$ is a {\em locally convex pleated plane}. More precisely, let $\cM\cL(\htwo)$ be measured geodesic laminations on $\htwo$ and $\cM\cL_0(\htwo) \subseteq \cM\cL(\htwo)$ the subspace of laminations with finite support. That is $\lambda \in \cM\cL_0(\htwo)$ if it is the union of a finite collection of disjoint geodesics $\ell_i$ with positive weights $\theta_i$. Then $\lambda$ determines a continuous map $p_\lambda\colon\htwo\to\hthree$, unique up to post-composition with isometries of $\hthree$, that is an isometry on the complement of the support of $\lambda$ and is ``bent'' with angle $\theta_i$ at $\ell_i$. By continuity we can extend this construction to a general $\lambda \in\cM\cL(\htwo)$. An exposition of the following theorem of Thurston can be found in \cite{KT-projective}.

\begin{theorem}\label{thurston_param}
Given $f \in P(\Delta)$ there exists maps $c_f\colon\Delta\to\htwo$ and $p_f\colon\htwo\to\hthree$ and a lamination $\lambda_f$ such that $p_f$ is a locally, convex pleated surface pleated along $\lambda_f$, $\Ep_{\rho_f} = p_f \circ c_f$ and the map $f\mapsto \lambda_f$ is a homeomorphism from $P(\Delta)\to\cM\cL(\htwo)$. Furthermore the maps $c_f\colon(\Delta, \rho_f)\to \htwo$ and $\Ep_{\rho_f}\colon(\Delta,\rho_f)\to \hthree$ are 1-Lipschitz.
\end{theorem}

\subsection{Average Bending Bound}
Average bending was introduced by the first author in the study of convex hulls of quasifuchsian groups (see \cite{Bridgeman:bending} and \cite{ BC:bending1}). This had  applications in the work of Epstein, Marden and Markovic in their paper   \cite{EMM1}. The idea of average bending is to relate the injectivity radius of the convex hull to the amount of bending per unit length along geodesic arcs. In their work, Epstein, Marden and Markovic, used an equivalent formulation of average bending,  called  {\em roundedness}. 
 
 Given $\lambda \in \mathcal{ML}(\Delta)$ and $\alpha$ a transverse arc, we let $\lambda(\alpha)$ be the $\lambda$-measure of $\alpha$. We then define the {\em average bending norm} to be
$$||\lambda||_L = \sup \{ \lambda(\alpha)\ | \alpha \mbox{ an open geodesic arc of length } L\}.$$

If $\lambda$ is a lift of a measured lamination on a closed hyperbolic surface, then $||\lambda||_L$ is bounded but in general $||\mu||_L$ may be infinite. For simplicity, we will let $\|\mu\|_1 = \|\mu\|$.

We have the following compactness result;

\begin{lemma}
Given $L, M > 0$ then the set $C(L,M) = \{ \lambda\ |\ ||\lambda||_L \leq M\}$ is precompact.
\end{lemma}

{\bf Proof:}
Let $G(\Delta)$ be the space of (unoriented) geodesics in the hyperbolic plane. We define the space of {\em geodesic currents} ${\mathcal C}(\Delta)$ to be the space of non-negative Borel measures on $G(\Delta)$ with the weak$^*$ topology. The topology on $\mathcal{ML}(\Delta)$ is  that of a closed subspace of ${\mathcal C}(\Delta)$. Given an open geodesic arc $\alpha$,  we let $U_\alpha \subseteq G(\Delta)$ be the set of all geodesics transverse to $\alpha$. We define 
$$\mathcal U = \{ U_\alpha\ | \ \alpha \mbox{ an open geodesic arc of length } L\}.$$

Then $\mathcal U$ is an open cover of $G(\Delta)$. 

We let $\mathcal K$ be the set of continuous functions on $G(\Delta)$ with  support  subordinate to the cover $\mathcal U$. Then  for each $\phi \in \mathcal K$ 
there exists a $U \in \mathcal U$ with $supp(\phi) \subset U$.  We have the map $I:\mathcal{ML}(\Delta) \rightarrow \R^{|\mathcal K|}$ given by $I(\lambda) = (\lambda(\phi))_{\phi\in \mathcal K}$. This map is a homeomorphism onto its image. 

If $\phi \in \mathcal K$ then there is a $U \in \mathcal U$ with  $supp(\phi) \subset U \in \mathcal U$. Therefore for  $\lambda \in C(M,L)$  
$$\lambda(\phi) \leq \lambda(U) \leq M.$$

Therefore  $C(M,L)$ is homeomorphic to a subset of  $[0,M]^{|\mathcal K|}$ which is compact by Tychanoff's theorem. Therefore $C(M,L)$ is precompact.\eproof

\begin{corollary}
Given $L,M >0$ there exists an $R$ such that if $f$ is a locally univalent map with $\|\lambda_f\|_L \leq M$ then $$\|\phi_f\|_\infty < R$$.
\label{compactness-arg}
\end{corollary}

{\bf Proof:}
We consider the family $F(M,L)$ of $\phi_f = S(f) \in Q(\Delta)$  with $\lambda_f \in C(M,L)$. Then by Thurston, $F(M,L)$ is the image of $C(M,L)$ under a homeomorphism. Therefore $F(M,L)$ is precompact and has compact closure $K(M,L)$. Therefore there is an $R  > 0$ such that for all $\lambda_f \in C(M,L)$ then 
$$|\phi_f(0)| \leq R/4.$$
Therefore $||\phi_f(0)|| \leq R$ for all $\lambda_f \in C(M,L)$. As the set $K(M,L)$ is invariant under isometries of $\Hp$ it follows that $||\phi_f||_\infty \leq R$ for all $\lambda_f \in C(M,L)$. 
\eproof

\subsection{Convex Hull of Cone Manifold}
In this section $M$ will be a conformally compact hyperbolic cone-manifold with incompressible boundary and  all cone angles $\le 2\pi$. We let $\phi$ be the quadratic differential on the conformal boundary given by  uniformization. In \cite{bromberg}, the second author studied   the convex core boundary  of $M$. This is given by taking the Epstein surface for the projective metric which we denote by $S$. By \cite[Proposition 6.5]{bromberg} the surface $S$ is an embedded locally convex surface in $M$ bounding an end $\mathcal E$ of $M$ homeomorphic to $S\times[0,\infty)$. Also $\mathcal E$ does not contain any cone axes in its interior. The surface $S$ has intrinsic hyperbolic metric and has a bending  lamination $\beta_{\phi}$. We identify the universal cover $\tilde S$ with the hyperbolic disk $\Delta$ and obtain a lamination $\tilde\beta_{\phi}$. 

First some elementary lemmas about  balls in hyperbolic cone-manifolds.

\begin{lemma}
Let $S$ be the unit sphere in $\R^3$. Let $(\theta,z)$ be cylindrical coordinates on $S$ and for  $0< t\leq 2\pi$ define the spherical cone-surface
$$S_t = \{ (\theta,z) \in S\ |\ 0\leq \theta \leq t  \}/(0,z) \sim (t,z).$$
If $p_1,p_2,p_3 \in S_t$ then $d(p_i,p_j) \leq 2\pi/3$ for some $i,j, i\neq j$.
\label{conepacking}\end{lemma}

{\bf Proof:}
Assume not. We first take the case of $t = 2\pi$. Then $S_t = S$ the unit sphere. Then  letting $B(p,r)$ be an open disk of radius $r$ about $p \in S$, we have
$$p_2,p_3 \in \overline{B(p_1,2\pi/3)}^c = B(-p_1,\pi/3).$$
It follows that $d_S(p_2,p_3) \leq 2\pi/3$ giving our contradiction.

For $t < 2\pi$ we take a fundamental wedge domain $W_t$  for $S_t$ in $S$ above, and can assume  the $p_i$ are in the interior. Then by the spherical case two of the points have $d_S(p_i,p_j) \leq 2\pi/3$. As $d_{S_t}(p_i,p_j) \leq d_S(p_i,p_j)$ we obtain our contradiction.
\eproof

We have the following elementary  calculation on half-spaces in $\Hs$;

\begin{lemma}
Let $f:\R_+ \rightarrow \R_+$ be given by
$$f(R)  = \cosh^{-1}\left(\frac{2\cosh(R)}{\sqrt{1+ 3\cosh^2(R)}}\right).$$
Let  $H_1, H_2, H_3$ be half-spaces  in $\Hs$ such that $H_i\cap B(x,R)$ are disjoint. If each $H_i$ intersects $B(x,r)$ then  $r \geq f(R)$.
\label{ballpack}
\end{lemma}

{\bf Proof}
Let $r_i$ be the distance from $H_i$ to $x$ and let $D_i = H_i\cap \partial B(x,R)$ have spherical radius $\theta_i$.   Then we have $r_i \leq r$ and $\theta_i \geq  \theta$ where $\theta$ is the spherical radius of  $D = H\cap \partial B(x,R)$ where $H$ is a half-space a distance $r$ from $x$.  Therefore as each $D_i$ contains a disk or radius $\theta$, if the $D_i$ are disjoint, then there are 3 disks of radius $\theta$ which are disjoint.

We show that $\theta \leq \pi/3$. We let $S = \partial B(x,R)$ have the spherical metric given by angle subtended at $x$. If $\theta > \pi/3$ then the centers  $p_i$ of $D_i$ satisfy $d(p_i,p_j) > 2\pi/3, i\neq j$ contradicting Lemma \ref{conepacking}.

We have  a right-angled hyperbolic triangle with sides $r, R$ and angle $\theta$ between. Let $l$ be the length of the other side. Then solving we have
$$\sinh(l) =\sinh(R).\sin(\theta) \leq \frac{\sqrt{3}\sinh(R)}{2}.$$ 
and by the hyperbolic Pythagorean formula  
$$\cosh(r) = \frac{\cosh(R)}{\cosh(l)} = \frac{\cosh(R)}{\sqrt{1+\sinh^2(l)}} \geq \frac{\cosh(R)}{\sqrt{1+ \frac{3}{4}\sinh^2(R)}} = \frac{2\cosh(R)}{\sqrt{1+ 3\cosh^2(R)}}.$$
\eproof

We now consider balls in our cone manifold $M$. We let  $\tilde M$ be the universal cover with convex hull $C(\tilde M)$. The end $\mathcal E$ lifts to $\tilde {\mathcal E}$ a component of the complement of $C(\tilde M)$ with boundary $\tilde S$. As $M$ has incompressible boundary, then $\pi_1(\tilde{\mathcal E})$ is trivial. 

The space $\tilde M$ is a hyperbolic cone manifold and the cone axes $\mathcal C$ lift to $\tilde{\mathcal C}$. For $p \in \tilde M$ we define balls in the usual way, i.e. $B(p,r) = \{ q \in \tilde M\ |\ d(p,q) \leq r\}$. We note that $B(p,r)$ may not be topologically a ball or isometric to a hyperbolic ball. For a point $p$, we define $r(p)$ to be the maximum radius such that $B(x,r(p))$ is embedded and isometric to a hyperbolic ball of radius $r(p)$. Note for $p\in \tilde{\mathcal C}$, $r(p) =0$ and otherwise $r(p) > 0$ and $r(p)$ equals is the injectivity radius of $p$ in $\tilde M-\tilde{\mathcal C}$. For $p \in \tilde M$ we further define $d(p)$ to be the minimum distance to the cone axes $\tilde{\mathcal C}$.

We first bound the average bending for points with $r(p)$ bounded below.

\begin{lemma}
Let $M$  be a conformally compact hyperbolic cone-manifold with incompressible boundary and  all cone angles $\le 2\pi$. Let $p \in \tilde S$ and $\alpha$ a closed geodesic arc on $\tilde S$ with midpoint $p$ and length less than $2f(r(p))$. Then 
$$\tilde\beta_{\phi}(\alpha) < 2\pi$$
\label{ballarg}
\end{lemma}

{\bf Proof:}
We let  $H_s$ be the 1-parameter family of support half-spaces from $\alpha(0)$ to $\alpha(1)$. We consider  $S = \partial B(p,r(p))$  and disks $D_s = H_s\cap S$.  We let $s_1$ be the smallest $s$ such that $D_0, D_s$ have disjoint interiors. Then we have $\tilde\beta(\alpha([0,s_1])) < \pi$. If there is no such $s_1$ then we have $\tilde\beta(\alpha) < \pi$ and we're done. 

We now let $s_2$ be the smallest $t$ such that $D_{s_1}, D_s$ have disjoint interiors. Again it follows that $\tilde\beta(\alpha[s_1,s_2]) < \pi$ giving $\tilde\beta(\alpha([0,s_2])) < 2\pi$. If no such $t_2$ exists then $\tilde\beta(\alpha) < 2\pi$ and we are also done. 

We first show that $D_0, D_{s_2}$ do not intersect. If $D_0, D_{s_2}$ do intersect, we extend $\alpha([0,s_2])$ to a closed curve $\alpha'$ by joining $\alpha(0), \alpha(s_2)$ by a piecewise geodesics on $\partial H_0\cup \partial H_{s_2}$.  We note that $\tilde{\mathcal E}$ is simply connected.  We get our contradiction by showing that curve $\alpha'$ in  $\tilde{\mathcal E}$  is homotopically non-trivial. 
The curve  $\alpha'$ is homotopic to a simple closed curve $\alpha''$ in $D_0\cup D_{s_1}\cup D_{s_2} \subset \partial B$ via a homotopy in $B\cap (H_0\cup H_{s_1}\cup H_{s_2}) \subset \tilde{\mathcal E}$. But as arc $\alpha$ is  transverse to a bending line $b$ then $\alpha''$ separates the points $b\cap \partial B$ in $\partial B$. Therefore $\alpha''$ is non-trivial in $\tilde M -b$. As $M$ has incompressible boundary, $\alpha''$ is trivial $\tilde M - C(\tilde M) \subset \tilde M- b$ and we obtain our contradiction. Thus if $s_2$ exists, then $D_0,D_{s_2}$ do not intersect. We then obtain a contradiction from the above lemma as $H_0,H_{s_1}, H_{s_2}$ are disjoint in $B(p,r(p))$ and intersect $B(p, f(r(p))$.
\eproof

We use the same argument as above to bound average bending for points close to the cone axes.

\begin{lemma}
Let $M$  be a conformally compact hyperbolic cone-manifold with incompressible boundary and  all cone angles $\le 2\pi$. Let $\tilde c \in \tilde{\mathcal C}$ have an embedded tube $U_{\tilde c}$ of radius $R$ and $p \in U_{\tilde c}$ with $d(p,c) \leq f(R)/2$. If  $\alpha$ is a closed geodesic arc on $\tilde S$ with midpoint $p$ and length less than $f(R)$ then 
$$\tilde\beta_{\phi}(\alpha) < 2\pi$$
\label{ballarg2}
\end{lemma}
{\bf Proof:}
We  let  $q \in \tilde{\mathcal C}$  be the nearest point of to $p$ on $\tilde{\mathcal C}$ and consider $B = B(q,R_0)$. If $\alpha$ is a geodesic arc of length $f(R_0)$ centered about $p$, then $\alpha$ is in $B(q,f(R_0))$. We let $S = \partial B$, then $S$ is a sphere with two cone points.  We again consider  $H_s$  the 1-parameter family of support planes from $\alpha(0)$ to $\alpha(1)$ and let $D_s = H_s \cap S$. Then $D_s$ are  disks in $S$ whose interior are disjoint from the cone points. Then analysing as in Lemma \ref{ballpack}, we obtain 3 disks with disjoint interiors on $S$. By Lemma \ref{conepacking} the disks cannot be disjoint which gives a contradiction. Thus we have $\tilde \beta(\alpha) < 2\pi$. 
\eproof

To bound our average bending uniformly for a given length, reduces now to showing that $r(p)$ is bounded away from zero for points far from the cone axes. This is the purpose of the following two lemmas.

\begin{lemma}
Let $M$  be a conformally compact hyperbolic cone-manifold with incompressible boundary and  all cone angles $\le 2\pi$. Let $\mathcal C$ be the cone-axes and for $\tilde c\in \tilde{\mathcal C}$ let $U_{\tilde c}$ be the $R$ neighborhood of $\tilde c$ in $\tilde M$. Let  $R$ be such that $U_{\tilde c}$ are embedded and disjoint. 
\begin{itemize}
\item If  $p \in \tilde M -  \cup_{\tilde c\in\tilde{\mathcal C}}  U_{\tilde c}$ then $r(p) \geq r(q)$ for some $q \in \partial  U_{\tilde c}$ and $\tilde c \in \tilde{\mathcal C}$. 
\item if $p \in  U_{\tilde c}$ then $d(p)= d(p,\tilde c)$   and 
\begin{eqnarray*}
r(p) =& d(p) &\theta_c \geq \pi,\\  
 \sinh(r(p)) =&   \sinh(d(p))\sin(\theta_c/2) & \theta_c \leq \pi.
 \end{eqnarray*}
 \end{itemize}
 
 \label{cat0}
\end{lemma}

{\bf Proof:}
We $\hat M$ to be the completion of the universal cover of $M-\mathcal C$. Then $\hat M$ is CAT(0). We let $\hat U_{\tilde c} \subseteq \hat M$ be the completions of the universal covers of $U_{\tilde c}$. 

Let $p \in \tilde M - \cup_{\tilde c\in \tilde{\mathcal C}}  U_{\tilde c}$. If $r(p)$ is achieved by an arc joining $p$ to an axis then $r(p) \geq R$ and as $r(q) \leq R$ for $q \in  \partial  U_{\tilde c}$ then the first statement follows. Otherwise there is a non-trivial geodesic $\gamma$ of length $2r(p)$ in $\tilde M -\tilde{\mathcal C}$. Therefore $\gamma$ must link a finite collection of  axes of elements $ \tilde{\mathcal C}$. If $\gamma$ links more than one axis then $\gamma$ is greater than the length of the shortest closed geodesic linking the axes. As this  is at least $2R$, then $r(p) \geq R$ as before. Therefore we can assume $\gamma$ links a single axis $\tilde c \in \tilde{\mathcal C}$. Then $\gamma$ lifts to $\hat\gamma$ a piecewise geodesic in $\hat M$ which is invariant under the action of the deck transformation $g_{\tilde c}$ corrseponding to $\tilde c$.
As  $\hat M$ is CAT(0) and the $\hat U_{\tilde c}$ are convex and complete, projection $\pi_{\tilde c}$ onto $\hat U_{\tilde c}$ is distance decreasing (see \cite[Proposition II.2.4]{Bridson:Haefliger:npc}). Thus as $\pi_{\tilde c}$ commutes with the action of $g_{\tilde c}$, the curve $\pi_{\tilde c}(\hat\gamma)$ descends to  a curve in $\tilde M$ of length $\leq 2r(p)$ which is contained in $  U_{\tilde c}$ linking $\tilde c$ with basepoint $q \in \partial U_{\tilde c}$. Thus $r(p) \geq r(q)$ and the first item is done.

For  $p \in \tilde U_{\tilde c}$ then trivially $d(p) = d(p,\tilde c)$.  We now describe the relation between $d(p),r(p)$. Then projecting as above, we have that $r(p)$ is attained by a curve in  $ U_{\tilde c}$. We take a fundamental domain for $ U_{\tilde c}$ to be a wedge $W$ of a hyperbolic tube of radius $R$ about a geodesic with wedge angle $t$ and have $p$ be on the central radial line of the wedge. Taking the largest  ball about $p$ that is embedded in $W$ it follows that  if $\theta_c \geq \pi$ then we have that $r(p) = d(p)$.   Otherwise $\theta_c \leq \pi$ and then $2r(p)$ is the length of the unique shortest geodesic arc with both endpoints $p$. Thus $r(p), d(p)$ are sides of a right angled triangle, with hypothenuse $d(p)$ and side of length $r(p)$ facing angle $\theta_c/2$. Therefore for $\theta_c \leq \pi$ by the hyperbolic sine formula
$$\sinh(r(p)) =   \sinh(d(p))\sin(\theta_c/2).$$
\eproof

\begin{lemma}
There is an explicit monotonic increasing function $g:\R_+\rightarrow \R_+$ such that the following holds. Let $M$ be a conformally compact hyperbolic cone-manifold with incompressible boundary and  all cone angles $\le 2\pi$ satisfying the conditions of Theorem \ref{infinity_bound}.
Then for $r \leq \frac{1}{2}\sinh^{-1}(\sqrt{2})$, if $p \in \partial C(\tilde M)$ with $d(p)\geq  r$  then $r(p) \geq g(r)$.
\label{bigballs}
\end{lemma}

{\bf Proof:}
By assumption for $c \in \mathcal C$, $c$ has embedded  tubular neighborhood $U_c$ of radius $R_c$ such that $R_c \geq \sinh^{-1}(\sqrt{2})$. We let   $R_0 = \sinh^{-1}(\sqrt{2})$. 
We lift the tubular neighborhoods to $\tilde M$ and  denote by $ U_{\tilde c}$  the lift for  $\tilde c \in \tilde{\mathcal C}$.
 
By the  Lemma \ref{cat0}, we need only consider points in the neighborhoods $U_{\tilde c}$.
Thus we let $p \in  U_{\tilde c}$. Again by Lemma \ref{cat0} if $\theta_c \geq \pi$ then $r(p) \geq d(p)$ giving $r(p) \geq r$.
Similarly for $\pi/2 \leq \theta_c \leq \pi$, we have 
$$\sinh(r(p)) = \sin(\theta_c/2)\sinh(d(p)) \geq \frac{1}{\sqrt{2}}\sinh(r).$$
This gives a bound on $r(p)$ for $\theta_c \geq \pi/2$.  

We now consider $\theta_c \leq \pi/2$. By \cite[Lemma 3.3]{bromberg} all support half-spaces are embedded in $\mathcal E$.  Let $H$ be a half space intersecting $ U_{\tilde c}$ with distance $d$ from the cone axis. We take a wedge fundamental domain with the nearest point of $H$ being centered. Then in order for $H$ to be embedded in $ U_{\tilde c}$, it cannot intersect the radial sides of the wedge. Therefore we must have $d> d_c$ where $d_c, R_c$ form a right-angled triangle with hypothenuse $R_c$ and angle between the sides $\theta_c/2$. Labeling the other side of the triangle $l$ we have by hyperbolic geometry (see \cite[formulas III.5, III.6]{FN})
$$\sinh(l) = \sinh(R_c)\sin(\theta_c/2)\qquad \tanh(l) = \sinh(d_c)\tan(\theta_c/2).$$
Thus for $p \in  U_{\tilde c}$ we have $d(p) \geq d_c$.  Therefore substituting
$$\sinh(r(p)) \geq \sinh(d_c)\sin(\theta_c/2) = \tanh(l)\cos(\theta_c/2) \geq \frac{1}{\sqrt{2}}\frac{\sinh(R_c)\sin(\theta_c/2)}{\sqrt{1+\sinh^2(R_c)\sin^2(\theta_c/2)}}$$
To obtain a  bound, we use our assumptions  in Theorem \ref{infinity_bound} and applying Theorem \ref{cone_margulis} we have 
$$\theta_c L_c\sinh(2R_c) = 1 \qquad\mbox{and}\qquad L_c \leq L_0\theta_c.$$
Therefore
$$\theta_c \geq \frac{1}{\sqrt{\sinh(2R_c)L_0}}.$$ 
It follows that 
$$\sinh(R_c)\sin(\theta_c/2)\geq \frac{1}{\sqrt{L_0}}  \frac{\sinh(R_c)}{\sqrt{\sinh(2R_c)}} \frac{\sin(\theta_c/2)}{\theta_c} \geq  \frac{1}{\sqrt{L_0}} 
 \sqrt{\frac{\tanh(R_c)}{2}}\frac{\sqrt{2}}{\pi}\geq \sqrt{\frac{\tanh(R_0)}{\pi^2  L_0}}.$$
As $\tanh(R_0) =\sqrt{2/3}$ and $x/\sqrt{1+x^2}$ is monotonic, then
$$\sinh(r(p)) \geq \frac{1}{\sqrt{2}}\frac{\sqrt{\frac{\tanh(R_0)}{\pi^2  L_0}}}{\sqrt{1+\frac{\tanh(R_0)}{\pi^2  L_0}}} =  \frac{1}{\sqrt{2+2\pi^2\cotanh(R_0) L_0}} \geq \frac{1}{\sqrt{2+24 L_0}}.$$
Thus for $p \in U_{\tilde c}$ we have  $r(p) \geq g(r)$ with 
$$\sinh g(r) =\min\left(\frac{1}{\sqrt{2+24L_0}},\frac{1}{\sqrt{2}}\sinh(r)\right)$$
 Thus  combining the  bounds, we have  $r(p) \geq g(r)$ with 
$$\sinh g(r) = \min\left(\frac{1}{\sqrt{2+24L_0}},\frac{1}{\sqrt{2}}\sinh(r)\right)$$
  \eproof

In \cite{BC:bending1}, the first author and Canary  proved the following.

\begin{theorem}[{\cite[Theorem 3]{BC:bending1}}]
Let $f:\Hp \rightarrow \Hs$ be an embedded convex pleated plane then its bending lamination $\beta_f$ satisfies
$$\|\beta_f\|_{L} < 2\pi$$
for $L \leq 2\sinh^{-1}(1).$
\label{bcbend}\end{theorem}

We now use the  Lemmas \ref{ballarg2} and \ref{bigballs} above to generalize Theorem \ref{bcbend}  for cone-deformations. 

\begin{prop}
Let $M$ be a conformally compact hyperbolic cone-manifold with incompressible boundary and  all cone angles $\le 2\pi$ satisfying the conditions of Theorem \ref{infinity_bound}. Then
$$||\beta_{\phi}||_{L} < 2\pi$$
for any  $L \leq 2f(g(f(\sinh^{-1}(\sqrt{2}))/2)) = .1529$.
\label{bend_explicit}
\end{prop}

{\bf Proof:}
 We let $r = f(R_0)/2$ where $R_0 = \sinh^{-1}(\sqrt{2})$. 

If $d(p) \geq  r$. Then as $f(x) \leq x$ we have
$$r = \frac{f(R_0)}{2} \leq \frac{R_0}{2} = \frac{1}{2}\sinh^{-1}(\sqrt{2}).$$
Therefore we  can apply   Lemma \ref{bigballs} to $p$ to get  $r(p) \geq g(r)$. Therefore for   $L =2f(g(r))$ then $\tilde\beta(\alpha) \leq 2\pi$ for any geodesic arc $\alpha$ of length less than $L$ centered at $p$. 

 If $d(p) \leq r =f(R_0)/2$, as $c\in \mathcal C$ has an embedded tubes of radius $R_c > R_0$ then by Lemma \ref{ballarg2}   if $\alpha$ is an arc of length $L \leq f(R_0)$ then  $\tilde\beta(\alpha) \leq 2\pi$.

Combining the bounds we have
$$\|\beta\|_L < 2\pi$$
For  $L \leq  \min(2r, 2f(g(r)))$. As $R_0 =\sinh^{-1}(\sqrt{2})$ and we can assume $L_0 < 1$, then  
$$\min(2r, 2f(g(r))) = 2f(g(r)) = .152958.$$ 
 \eproof

We now prove the main result of this section.

{\bf Proof of Theorem \ref{infinity_bound}:} 
By the above, there exists an $L$ such that 
  $$\|\beta_{\phi}\|_{L} < 2\pi.$$
Therefore by Corollary \ref{compactness-arg} 
$$\|\Sigma\|_\infty \leq K$$
for some $K$ universal. 
\eproof

\section{Proof of Theorems \ref{inf_cone} and  \ref{cone}}
We now bring our work together to prove the main results of the paper. Before doing so we will need to summarize the necessary results about deformations of cone-manifolds. As in the introduction we have a compact 3-manifold $\bar N$ with a collection $\mathcal C$ of disjoint, simple closed curves in the interior. We will examine a family of conformally compact hyperbolic cone-manifold structures on $\bar N$ with cone locus $\mathcal C$.

\begin{theorem}[\cite{bromberg}]\label{cone_estimates}
Let $M_t$ be a one parameter family of cone-manifolds given by Theorem \ref{cone_deformation} and let $L_c(t)$  be the length of a component $c$ of $\mathcal C$ in $M_t$ and $L_{\mathcal C}$ the sum of the $L_c$.
\begin{itemize}
 \item $$L_c(t) \le \frac{tL_c(2\pi)}{\pi}$$

\item In each $M_t$ there is a union $U_t$ of embedded, disjoint tubular neighborhoods of the components of $\mathcal C$ of radius $\ge \sinh^{-1}\sqrt 2$.

\item The time $t$ derivative of $M_t$ is represented by a model Hodge form $\omega_t$ with
$$\int_{M_t\backslash U_t} ||\omega_t||^2 \leq \frac{3}{14} \cdot \frac{L_{\mathcal C}(t)}{t}\le\frac{3L_{\mathcal C}(2\pi)}{14\pi}.$$
\end{itemize}
\end{theorem}

Note that the statement in the final bullet is not the actual statement of Proposition 4.2 in \cite{bromberg} but rather a direct application of the first inequality of the proof where we assume that the radius of the tubular neighborhoods is $\sinh^{-1}\sqrt 2$ rather than the larger radii assumed in that proposition.

We are now ready to prove our main theorem bounding  the $L^2$-norm of the derivative of the path of complex projective structures.

{\bf Proof of Theorem \ref{inf_cone}:} We assume $t$ has been fixed throughout the proof.

For the path $\Sigma_t$ of complex projective structures on the boundary of $M_t$, by Theorem \ref{infinity_bound} we have that $\|\Sigma_t\|_\infty \le K$. Therefore by Theorem \ref{epstein_hyp} there is a convex surface $S$ in $M_t$ cutting of an end $\E$ such that $(1+2K)\hat g_X$ is the metric at infinity for $S$. Note that while $\E$ will be disjoint from the cone locus in $M_t$ it may intersect the tubular neighborhood $U_t$ of the cone locus. To correct this we need to remove the collar of width $\sinh^{-1}\sqrt 2$ from $\E$. This is the end $\E_\eta$ where $\eta = e^{-\sinh^{-1}\sqrt 2}$. 

By Theorem \ref{cone_estimates} we have that
$$ \int_{M_t\backslash U_t} ||\omega_t||^2 \le\frac{3L_{\mathcal C}}{14\pi}$$
and since $\E_\eta \subset M_t\backslash U_t$ this implies that
$$\int_{\E_\eta} \|\omega_t\|^2  \le\frac{3L_{\mathcal C}}{14\pi}.$$
As $\omega_t$ is a model Hodge form Theorem \ref{cone_estimates} implies that
$$\|\Phi_t\|^2_{(1+2K)\hat g_X,2} \le \frac{1}{8\eta^2} \frac{3L_{\mathcal C}}{14\pi}.$$
As
$$\|\hat \Phi_t\|^2_{\hat g_X,2}=(1+2K) \|\Phi_t\|^2_{(1+2K)\hat g_X,2} $$
this gives 
$$\|\Phi_t\|_{\hat g_X,2} \le c_{\rm drill} \sqrt{L_\mathcal C}$$
where
$$c_{\rm drill} = \frac{1}{4\eta} \sqrt{\frac{3(1+2K)}{7\pi}}.$$
\eproof

Our main results now follows immediately.

{\bf Proof of Theorem \ref{cone}:} Integrating the above, we get the $L^2$-bound
$$||\Phi(\Sigma_0,\Sigma_{2\pi})||_2  \leq \int_{0}^{2\pi} ||\Phi_t||_2 dt \leq   2\pi c_{\rm drill}\sqrt{L_{\mathcal C}}.$$
\eproof

\bibliography{bib,math}

\def\cprime{$'$}
\begin{thebibliography}{EMM}

\bibitem[Bri]{Bridgeman:bending}
M.~Bridgeman.
\newblock {Average bending of convex pleated planes in hyperbolic three-space}.
\newblock {\em Invent. Math.} {\bf 132}(1998), 381--391.

\bibitem[BBB]{wp-paper}
M.~Bridgeman, J.~Brock, and K.~Bromberg.
\newblock {The Weil-Petersson gradient flow of renormalized volume and
  3-dimensional convex cores}.
\newblock preprint (2020), https://arxiv.org/abs/2003.00337.

\bibitem[BC1]{BC:bending1}
M.~Bridgeman and R.~Canary.
\newblock {Bounding the bending of a hyperbolic three-manifold}.
\newblock {\em Geom. Dedicata} {\bf 96}(2003), 211--240.

\bibitem[BW]{bridgeman2019uniform}
M.~Bridgeman and Y.~Wu.
\newblock {Uniform bounds on harmonic {B}eltrami differentials and
  {W}eil-{P}etersson curvatures}.
\newblock {\em J. Reine Angew. Math.} {\bf 770}(2021), 159--181.

\bibitem[BH]{Bridson:Haefliger:npc}
M.~Bridson and A.~Haefliger.
\newblock {\em Metric Spaces of Non-Positive Curvature}.
\newblock Springer-Verlag, 1999.

\bibitem[BB]{Brock:Bromberg:inflexible}
J.~Brock and K.~Bromberg.
\newblock {Geometric inflexibility and 3-manifolds that fiber over the circle}.
\newblock {\em Journal of Topology} {\bf 4}(2011), 1--38.

\bibitem[Bro]{bromberg}
K.~Bromberg.
\newblock {Hyperbolic cone-manifolds, short geodesics and Schwarzian
  derivatives}.
\newblock {\em J. Amer. Math. Soc.} {\bf 17}(2004), 783--826.

\bibitem[Brm]{bromberg:thesis}
K.~Bromberg.
\newblock {Rigidity of geometrically finite hyperbolic cone-manifolds}.
\newblock {\em Geom. Dedicata} {\bf 105}(2004), 143--170.

\bibitem[Cal]{calabi:rigidity}
Eugenio Calabi.
\newblock {On compact, {R}iemannian manifolds with constant curvature. {I}}.
\newblock In {\em Proc. {S}ympos. {P}ure {M}ath., {V}ol. {III}}, pages
  155--180. American Mathematical Society, Providence, R.I., 1961.

\bibitem[Eps]{epstein-envelopes}
C.~Epstein.
\newblock {Envelopes of horospheres and Weingarten surfaces in hyperbolic
  3-spaces}.
\newblock preprint Princeton University, (1984), available at
  \verb+https://www.math.upenn.edu/~cle/papers/WeingartenSurfaces.pdf+.

\bibitem[EMM]{EMM1}
D.~B.~A. Epstein, A.~Marden, and V.~Markovic.
\newblock {Quasiconformal homeomorphisms and the convex hull boundary}.
\newblock {\em Ann. of Math. (2)} {\bf 159}(2004), 305--336.

\bibitem[FN]{FN}
Werner Fenchel and Jakob Nielsen.
\newblock {\em Discontinuous Groups of Isometries in the Hyperbolic Plane}.
\newblock De Gruyter, 2011.

\bibitem[HK]{Hodgson:Kerckhoff:rigidity}
C.~Hodgson and S.~Kerckhoff.
\newblock {Rigidity of hyperbolic cone-manifolds and hyperbolic Dehn surgery}.
\newblock {\em J. Diff. Geom.} {\bf 48}(1998), 1--59.

\bibitem[KT]{KT-projective}
Yoshinobu Kamishima and Ser~P. Tan.
\newblock {Deformation spaces on geometric structures}.
\newblock In {\em Aspects of low-dimensional manifolds}, volume~20 of {\em Adv.
  Stud. Pure Math.}, pages 263--299. Kinokuniya, Tokyo, 1992.

\bibitem[KS]{KS08}
K.~Krasnov and J-M. Schlenker.
\newblock {On the renomalized volume of hyperbolic \hbox{3-manifolds}}.
\newblock {\em Comm. Math. Phys.} {\bf 279}(2008), 637--668.

\bibitem[Weil]{Weil:compact}
A.~Weil.
\newblock {On discrete subgroups of Lie groups}.
\newblock {\em Annals of Math.} {\bf 72}(1960), 369--384.

\end{thebibliography}

\bibliographystyle{math}
\end{document}